\documentclass[preprint,12pt]{elsarticle}

\usepackage{graphicx}
\usepackage{amsmath}
\usepackage{amssymb}
\usepackage{mathtools}
\usepackage{algorithm}
\usepackage[noend]{algpseudocode}
\usepackage{scalerel}
\usepackage{url}
\usepackage{lineno}

\newcommand\equref[1]{Eq.~\ref{eq:#1}}
\newcommand\figref[1]{Fig.~\ref{fig:#1}}
\newcommand\algoref[1]{Algorithm~\ref{alg:#1}}

\DeclarePairedDelimiter{\rbra}{(}{)}
\DeclarePairedDelimiter{\cbra}{\{}{\}}
\DeclarePairedDelimiter{\abs}{|}{|}
\DeclarePairedDelimiterX{\Set}[2]{\lbrace}{\rbrace}{#1\,\delimsize\vert\,#2}

\DeclareMathOperator{\sech}{sech}

\journal{{}}

\begin{document}

\begin{frontmatter}

\title{
Fast parallel calculation of modified Bessel function\\
of the second kind and its derivatives
%of the second order
}

\author{Takashi Takekawa}
\ead{takekawa@cc.kogakuin.ac.jp}
\address{
Faculty of Informatics, Kogakuin University of Technology and Engineering,\\
1-24-2 Nishi-Shinjuku, Shinjuku-ku, Tokyo 163-8677, JAPAN
}

\begin{abstract}
There are three main types of numerical computations
for the Bessel function of the second kind:
series expansion, continued fraction, and asymptotic expansion.
In addition, they are combined in the appropriate domain for each.
However, there are some regions where the combination of these types requires sufficient computation time to achieve sufficient accuracy,
however, efficiency is significantly reduced when parallelized.
In the proposed method, we adopt a simple numerical
integration concept of integral representation.
We coarsely refine the integration range beforehand,
and stabilize the computation time by performing the integration 
calculation at a fixed number of intervals. Experiments demonstrate 
that the proposed method can achieve the same level of accuracy
as existing methods in less than half the computation time.
\end{abstract}

\begin{keyword}
Bessel functions
\sep
Numerical integration
\sep
Parallel execution
\end{keyword}

\end{frontmatter}

\section{Introduction}
The modified Bessel function of the second kind $K_v\rbra{x}$
is an important special function adopted in various fields \cite{watson}.
In addition to mathematics and physics,
it has become increasingly important in the fields of statistics and economics.
For example, normal inverse Gaussian
\cite{barndorff1997processes} has been garnering increasing attention
in economics \cite{barndorff1997normal}, biometrics \cite{hassan2017automated, das2016classification}, and machine learning \cite{o2016clustering, takekawa2020clustering}.

Three numerical methods exist calculating the Bessel function $K_v$: series expansion--based \cite{temme}, continued fraction-based \cite{temme, amos, cambell, thompson}, and asymptotic expansion-based methods \cite{olver, olver1}.
Because each computation method exhibits different accuracies
and computation time characteristics for $v$ and $x$, a combined implementation of these methods is generally adopted \cite{gsl, boost, tfp}. However, in some regions, an enormous amount of computation time is required for any method to achieve sufficient accuracy.

Recently, it has become common to adopt GPUs and multi-core processors to perform large-scale paralell or vectorized processing, and
it is desirable to support parallel processing for the computation of special functions.
When computing function values for several orders and variable combinations of variables simultaneously, computation time depends on the slowest computation combination. Therefore, if computation is slow for a particular combination, the overall performance will decrease significantly.

In this study, we propose an integration-based method for computing $K_v\rbra{x}$, which is suitable for parallelization.
Computing directly from the definition of the integral representation is simple and powerful \cite{schwartz}.
However, the method that involves determining the bin size and sustaining addition until convergence cannot handle a wide range of orders $v$ and arguments $x$. In particular, if the bin size is fixed, the number of additions until convergence varies 
significantly for each argument, making the method unsuitable for parallelization.

Therefore, in the proposed method, the range to be integrated is calculated beforehand, and the number of bins is fixed, instead of the width of the bin. Even if the evaluation of the integration range is coarse, it does not significantly affect the accuracy, and the integration can be performed
by adopting the parallelization feature.

The derivative of an arbitrary Bessel function can also be computed similarly. Because there are no major libraries
for the differentiation of Bessel functions of any order, it would be beneficial
to provide the implementations of these functions.

We implemented the proposed method using Python and TensorFlow \cite{tf}. The obtained results indicate that the proposed method outperforms conventional methods in terms of both accuracy and speed.

\section{Proposed methods}
For $v \in \mathbb{R}$, $x > 0$,
modified Bessel functions of the second kind can be represented
by the evaluation of integrals of the form \cite{watson}:
\begin{equation}
    K_v\rbra{x} = \int_0^\infty \!\! {f_{v, x}\rbra{t}}\;dt,
\end{equation}
where
\begin{gather}
    f_{v, x}\rbra{t} = \cosh\rbra{v t} \exp \rbra*{- x \cosh t}.
\end{gather}

\subsection{Shape of $f_{v,x}$}
Considering the shape of $f_{v,x}\rbra{t}$,
we examine the increase or decrease in its logarithmic form
\begin{gather}
    \log f_{v,x}\rbra{t} = g_{v,x}\rbra{t} = \log \cosh \rbra{v t} - x \cosh \rbra{t}.
\end{gather}
The derivatives of $g_{v,x}\rbra{t}$ relative to $t$ are given by
\begin{gather}
    g'_{v,x}\rbra{t} = v \tanh \rbra{v t} - x \sinh \rbra{t},\\
    g''_{v,x}\rbra{t} = v^2 \sech^2 \rbra{v t} - x \cosh \rbra{t},\\
    g'''_{v,x}\rbra{t} = - v^3 \sech^2 \rbra{vt} \tanh \rbra{vt} - x \sinh \rbra{t}.
\end{gather}
In particular,
$g_{v,x}\rbra{0} = -x$, $g'_{v,x}\rbra{0} = 0$ and $g'''_{v,x}\rbra{t} \le 0$ always hold.

\begin{figure}
    \centering
    \includegraphics{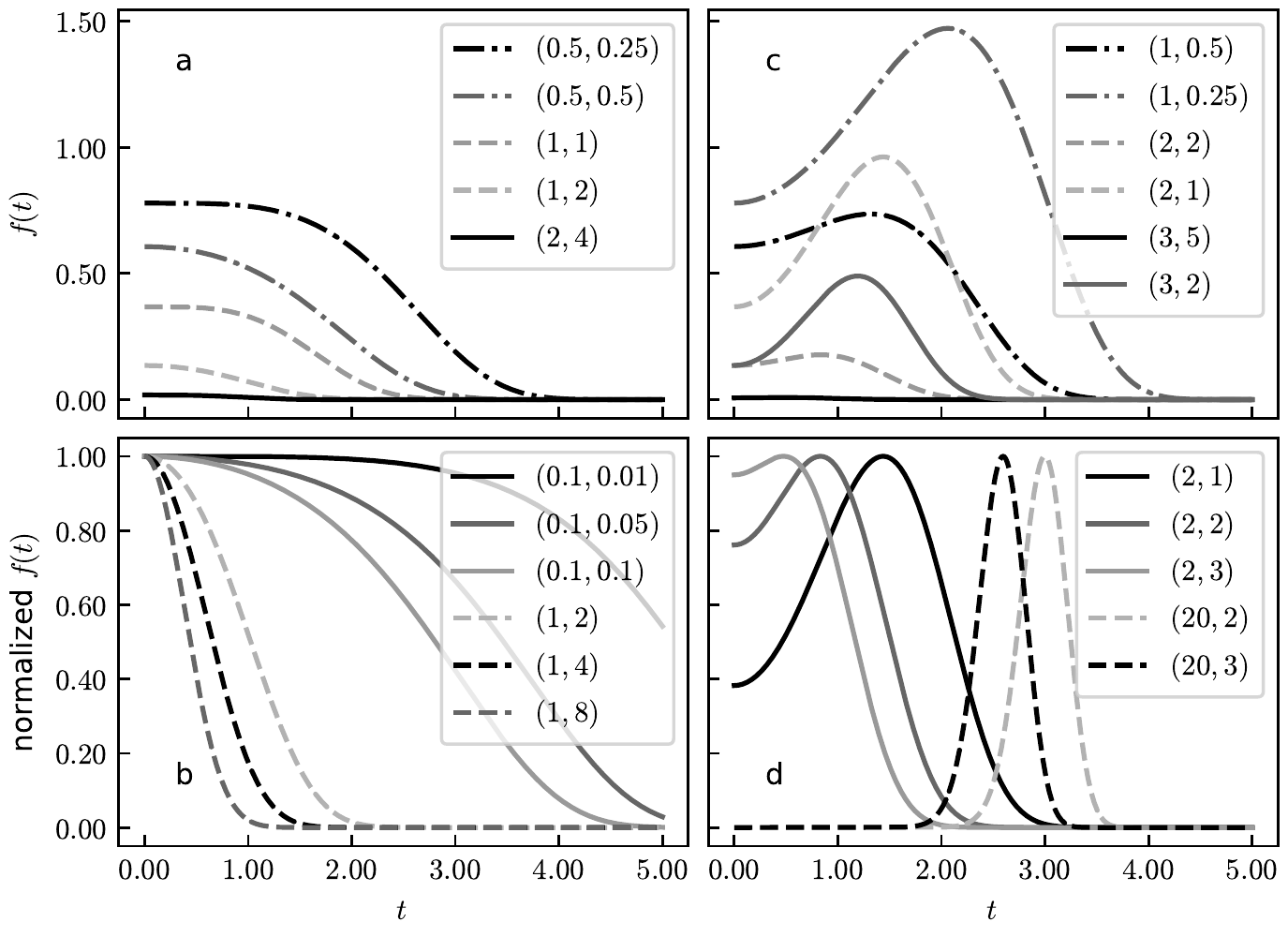}
    \caption{Shape of $f_{v,x}\rbra{x}$: two top panels (a, c) are $f_{v,x}\rbra{t}$ and two bottom panels (b, d) present the normalized form $f_{v,x}\rbra{t} / f_{v,x}\rbra{t}$. The two numbers in the legend are the order $v$ and argument $x$. The left and right panels (a, b) and (c, d) present examples with peaks at $t=0$, and $t>0$, respectively.}
    \label{fig:fvx}
\end{figure}

Therefore, in the case of $v^2 \le x$, 
$g_{v, x}\rbra{t}$ monotonically decreases
because $g''_{v,x}\rbra{0} \le 0$
(\figref{fvx}ab).
However, in the case of $v^2 > x$,
$g_{v, x}\rbra{t}$ has only one peak $t_p \ge 0$,
because $g''_{v,x}\rbra{0} > 0$
(\figref{fvx}cd).

When integrating $f_{v,x}\rbra{t}$ numerically,
we sololy need to consider the region where
$g_{v,x}\rbra{t} \ge \epsilon g_{v,x}\rbra{t_p}$,
in which $g_{v,x}\rbra{t}$ takes the maximum value at $t_p$
and $\epsilon$ denotes the machine epsilon.
In fact, from the shape of $g_{v,x}\rbra{t}$,
the region that satisfies the condition
can be defined as a single continuous range as:
\begin{align}
    [t_0, t_1]
    &= \Set*{t}{f_{v,x}\rbra{t} \ge \epsilon f_{v,x}\rbra{t_p}}\\
    &= \Set*{t}{g_{v,x}\rbra{t} \ge g_{v,x}\rbra{t_p} + \log \epsilon}.
\end{align}

\subsection{Find peak $t_p$}\label{secfr}
To determine the integral range,
we first find the maximum value.
In the case of $v^2 \le x$,
$t_p=0$ is obvious from the shape of $g_{v,x}\rbra{t}$.
In the case of $v^2 > x$,
$g_{v,x}\rbra{t}$ exibits the maximum at $t_p>0$.

To search for $t_p > 0$,
first, the range where the zero point $t_p$ of $g'_{v,x}\rbra{t_p}=0$ exists
is determined using the property that 
$g'_{v,x} > 0$ for $t<t_p$ and $g'_{v,x}\rbra{t} < 0$ for $t > t_p$.
Specifically,
for $m \in \mathbb{Z}$,
 to $[2^{m-1},2^m]$
find
the smallest $m$, such that $g'_{v,x}\rbra{2^m}<0$.
The details of this procedure (\texttt{FindRange}) are provided in \algoref{findrange}.

\begin{algorithm}
    \caption{Find range method}\label{alg:findrange}
    \begin{algorithmic}[1]
        \Require {$g\rbra{t_0} \ge 0$}
        \Function {FindRange}{$g, t_0$}
        \Comment{return $t_* \ge t_0$, such that $g\rbra{t_*} < 0$}
            \State $m \gets 0$
            \While {$g\rbra{t_0 + 2^{m+1}} \ge 0$}
                \State $m \gets m + 1$
            \EndWhile
            \State \Return $t_0 + 2^m$, $t_0 + 2^{m+1}$
        \EndFunction
    \end{algorithmic}
\end{algorithm}

Next, for the obtained existence range $[2^{m-1}, 2^m]$, find $t_p$, such that $f'\rbra{t_p} = 0$, using a combination of
the binary search and Newton methods. The details of this procedure (\texttt{FindZero}) are provided in \algoref{findzero}.

\begin{algorithm}
    \caption{Modified Newton Method}\label{alg:findzero}
    \begin{algorithmic}[1]
        \Require {$g(t_0)>0$, $g(t_1)<0$}
        \Function {FindZero}{$g, t_0, t_1$}
        \Comment{return $t_0 \le t_* \le t_1$, such that $g\rbra{t_*} = 0$}
            \State $t_0^{(0)}, t_1^{(0)}, i \gets t_0, t_1, 0$
            \Repeat
            \Comment{$\text{tol} = 1$ and $\text{max\_iter} = 10$}
            
                \State $t_\text{shrink} \gets t_0^{(i)} + \rbra{t_1^{(i)} - t_0^{(i)}}/2$\Comment{Binary search}
                \State $t_\text{newton} \gets \Call{Clip}{t_0^{i+1} - g\rbra{t_0^{(i)}}/g'\rbra{t_0^{(i)}}, t_0, t_\text{shrink}}$
                \Comment{Newton method}
                %\If {$t_0^{(i)} < t_\text{newton} < t_1^{(i)}$}
                %    \State $t_0^{(i+1)} \gets t_\text{newton}$
                %\EndIf
                \If {$g\rbra{t_\text{shrink}} < 0$}
                    \State $t_0^{(i+1)}, t_1^{(i+1)} \gets t_\text{shrink}, t_1^{(i)}$
                \ElsIf {$g\rbra{t_\text{newton}} < 0$}
                    \State $t_0^{(i+1)}, t_1^{(i+1)} \gets t_\text{newton}, t_\text{shrink}$
                \Else
                    \State $t_0^{(i+1)}, t_1^{(i+1)} \gets t_0^{(i)}, t_\text{newton}$
                \EndIf
                \State $i \gets i + 1$
            \Until {$\abs{t_0^{i} - t_0^{i-1}} < \text{tol}$ and $i < \text{max\_iter}$}
            \State \Return $t_0^{(i)}$
        \EndFunction
    \end{algorithmic}
\end{algorithm}

\subsection{Find the integration range $[t_0, t_1]$}\label{secfz}
In the cases of $v^2 \le x$
or of $v^2 > x$ and $f\rbra{0} - \epsilon f\rbra{t_p} > 0$,
$t_0=0$ from the shape of $f_{v,x}\rbra{t}$.
Otherwise, \texttt{FindZero} is applied
within the range $[t_0, t_p]$ to obtain
$t_0$, such that
$g_{v,x}\rbra{t_0} - g_{v,x}\rbra{t_p} - \log \epsilon = 0$.

For $t_1$, after applying \texttt{FindRange}
to determine the range $[t_p + 2^{m-1}, t_p + 2^m]$,
apply \texttt{FindZero}
within the range
to obtain
$t_1$, such that
$g_{v,x}\rbra{t_1} - g_{v,x}\rbra{t_p} - \log \epsilon = 0$.

\subsection{Integration}
After obtaining $t_0$ and $t_1$,
integrate numerically over the range with the fixed number of divisions $n$:
\begin{gather}
    K_v(x) \sim h \sum_{m=0}^{n} c_m f_{v, x}\rbra*{t_m},
\end{gather}
where
\begin{gather}
    h = \frac{t_1 - t_0}{n},\quad
    t_m = t_0 + mh, \\
    c_0 = c_n = 1/2, \quad
    c_m = 1 \,\, \rbra{m = 1, \ldots, n - 1}.
\end{gather}
In fact, for numerical stability, \texttt{log\_sum\_exp} was applied to $g_{v,x}\rbra{t_m}$:
\begin{gather}
    \log K_v(x) \sim 
    g_{v,x}\rbra{t_p} + \log \sum_{m=0}^{n} h \exp \cbra*{c_m \rbra*{g_{v, x}\rbra{t_m} - g_{v, x}\rbra{t_p}}}.\label{int}
\end{gather}

\subsection{Implementation}

The entire proposed algorithm based on integration,
which is hereafter denoted as ''I,''
is summarized in \algoref{main}.
Our implementations using Python and TensorFlow are available at \url{https://github.com/tk2lab/logbesselk}.
In section \ref{eval}, the proposed method is compared with implementations by series expansion (''S;'' see \ref{series}), continued fraction (''C;'' see \ref{continuedfraction}), and asymptotic expansion (''A;'' see \ref{asymptotic}).
The implementations adopted for the comparison are also available.
These implementations used for the comparison are also publicly available at the same already provided url.

\begin{algorithm}
    \caption{Main algorithm of logbesselk}\label{alg:main}
    \begin{algorithmic}[1]
        \Function {LogBesselK}{$v, x$}
        \Comment{return $\log K_v\rbra{x}$}
            \State $t_p \gets 0$
            \If {$g''_{v,x}\rbra{0} > 0$}\Comment{find $t_p > 0$}
                \State $t_s, t_e \gets \Call{FindRange}{g'_{v, x}\rbra{t}, 0}$
                \State $t_p \gets \Call{FindZero}{g'_{v, x}\rbra{t}, t_e, t_s}$
                %\Comment{$g_{v,x}\rbra{\delta} > 0$}
            \EndIf
            \State $t_0 \gets 0$
            \If {$g_{v,x}\rbra{0} - g_{v,x}\rbra{t_p} \le \log \epsilon$}\Comment{find $0 < t_0 \le t_p$}
                \State $t_0 \gets \Call{FindZero}{g_{v,x}\rbra{t} - g_{v,x}\rbra{t_p} - \log \epsilon, 0, t_p}$
            \EndIf
            \State $t_s, t_e \gets \Call{FindRange}{g_{v,x}\rbra{t}  - g_{v,x}\rbra{t_p} - \log \epsilon, t_p}$
            \State $t_1 \gets \Call{FindZEro}{g_{v,x}\rbra{t}  - g_{v,x}\rbra{t_p} - \log \epsilon, t_e, t_s}$
            \Comment{find $t_1 > t_p$}
            \State \Return $\int_{t_0}^{t_1} g_{v,x}\rbra{t} dt$\Comment{see \eqref{int}}
        \EndFunction
    \end{algorithmic}
\end{algorithm}

\section{Evaluation}\label{eval}

We adopted the calculations by Mathematica (Walfram) as an accurate reference
to evaluate accuracy in the range $v \in [0, 99]$, $x \in [10^{-1}, 10^{2.1}]$.
%$v = \Set*{10^\frac{i}{40} - 1}{i=0,1,\ldots,80}$,
%$x = \Set*{10^\frac{i}{40}}{i=0,1,\ldots,124}$
In this study, the errors were evaluated using $\log_{10} \rbra*{\abs{\Delta}/\epsilon + 1}$,
where $\Delta$ is the deviation from the reference.
The experiments are conducted in the following environment:
Intel(R) Xeon(R) W-2123 3.60GHz CPU,
128GB memory,
NVIDIA TITAN RTX GPU,
Ubuntu 20.04.2 LTS,
Python 3.8.10,
and TensorFlow 2.6.0.

\subsection{Accuracy}

\begin{figure}
    \centering
    \includegraphics{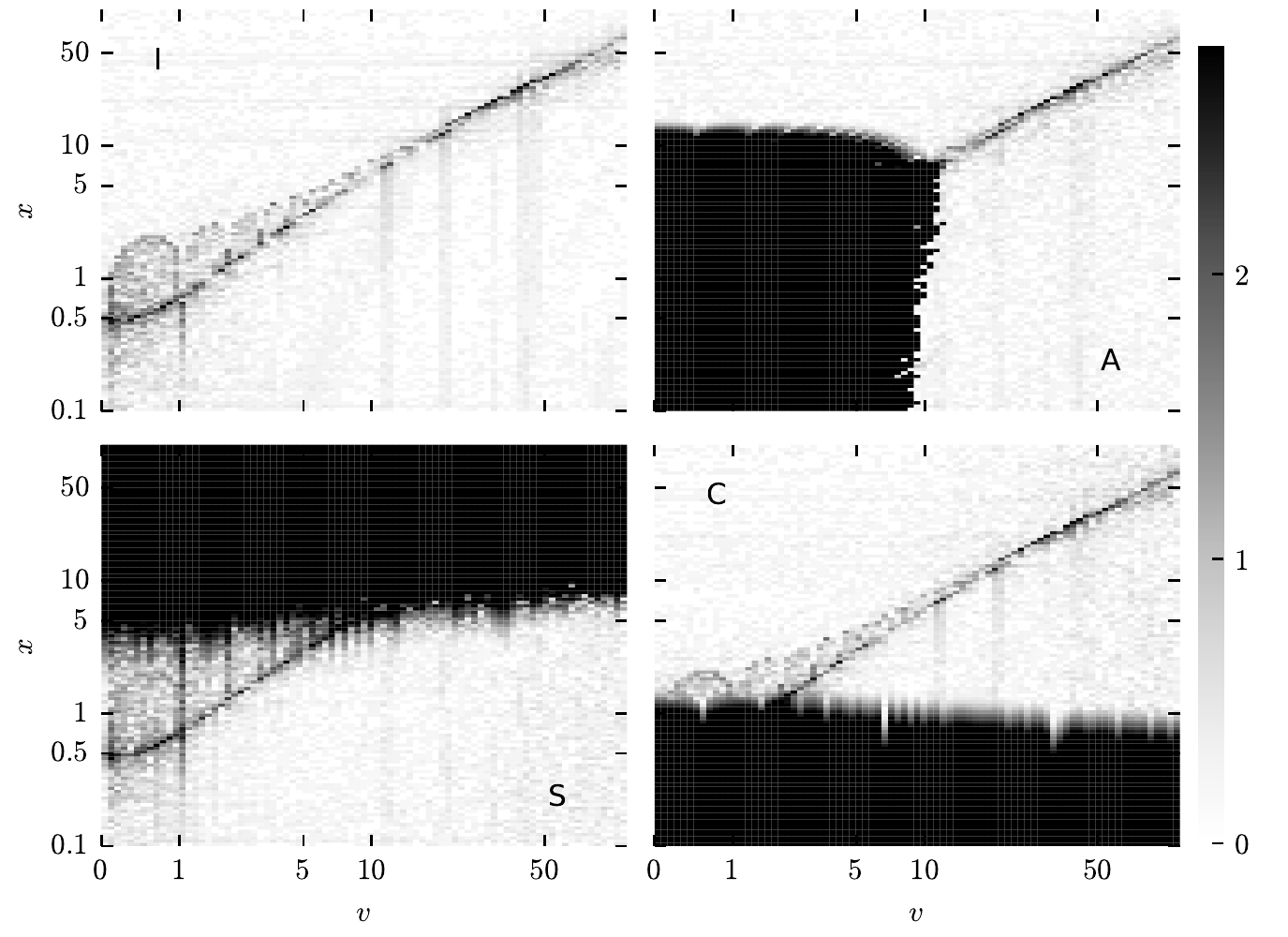}
    \caption{Error corresponding to order $v$ and argument $x$: I) the proposed integration method; S) series expansion method; C) continued fraction method; A) asymptotic expansion method. Error is defined by $\log \rbra{\abs{\Delta}/\epsilon+1}$, where $\Delta$ denotes the deviation from the reference. }
    \label{fig:acc}
\end{figure}

We first evaluate the accuracy of 
the methods based on series expansion (S),
continued fraction (C), asymptotic expansion (A),
and integration (I).
In the region near $v=x$, the accuracy of these methods deteriorates \cite{gil}.
It can be deduced that series expansion is accurate for $x < 5$, continued fraction is accurate for $x > 1$, and asymptotic expansion is accurate
in the range $v > 10$ or $x > 10$ (\figref{acc}).
However, there is no region where integration is inaccurate.

\subsection{Computational time}
\begin{figure}
    \centering
    \includegraphics{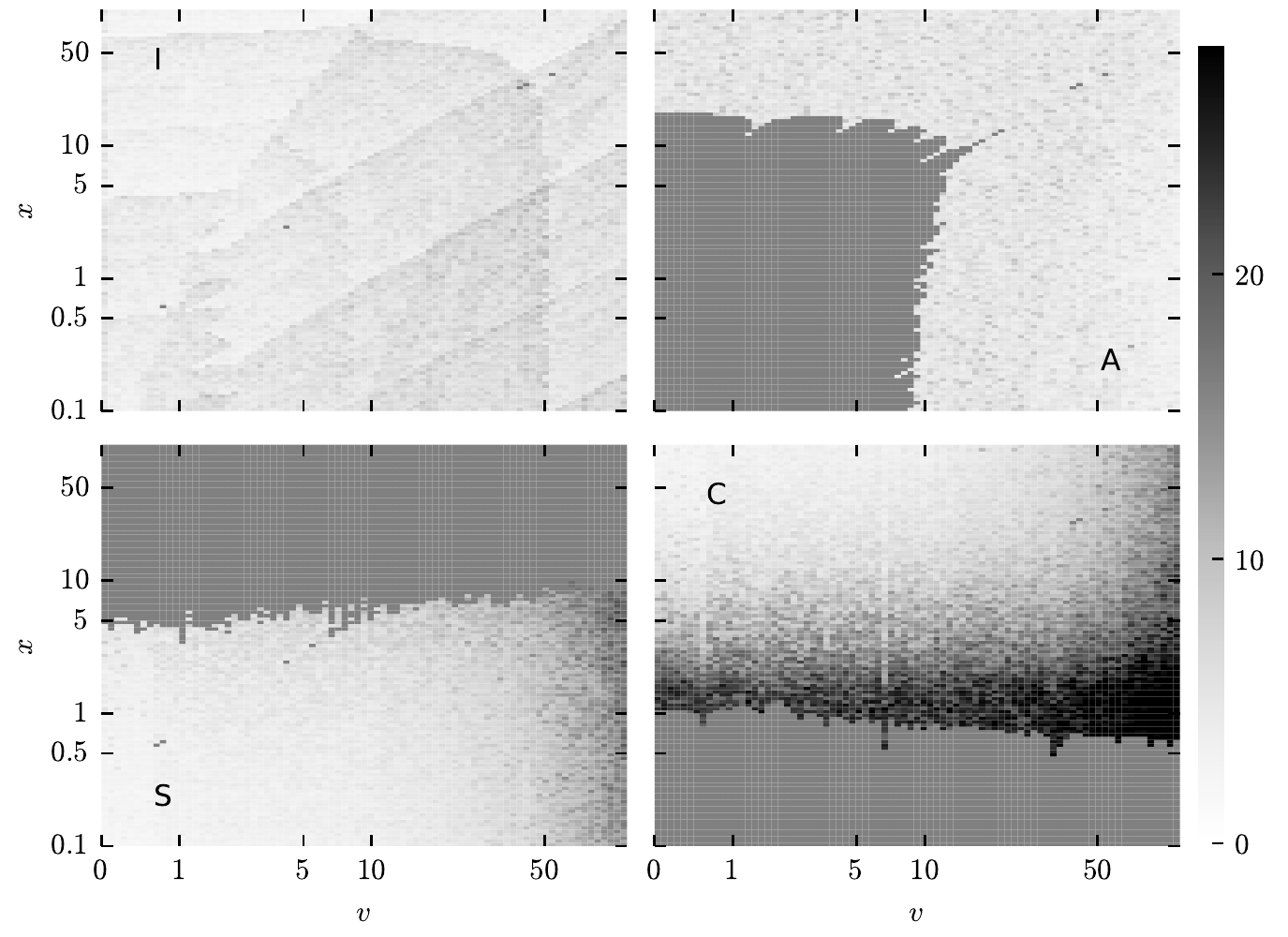}
    \caption{Computational time (msec) corresponding to order $v$ and argument $x$: I) the proposed integration method; S) series expansion method; C) continued fraction method; A) asymptotic expansion method. The computation time is not measured in the gray area because of the large error.}
    \label{fig:time}
\end{figure}

Next, we examined the computation time for each order $v$ and argument $x$ separately.
The computation time was obtained by excluding regions with errors of 4 or more.
 For I, the computation time tends to be generally uniform; however, for S, C, and A, there are regions in the domain where the computation time becomes large (\figref{time}).

 Although TensorFlow can compute multiple pairs of $v$ and $x$ simultaneously, the computation time in this case is the worst computation time for the included $v$ and $x$. Therefore, the worst-case computation time within the range assumed in the application is required to be small, instead of the average computation time.

\subsection{Combinations}

\begin{figure}
    \centering
    \includegraphics{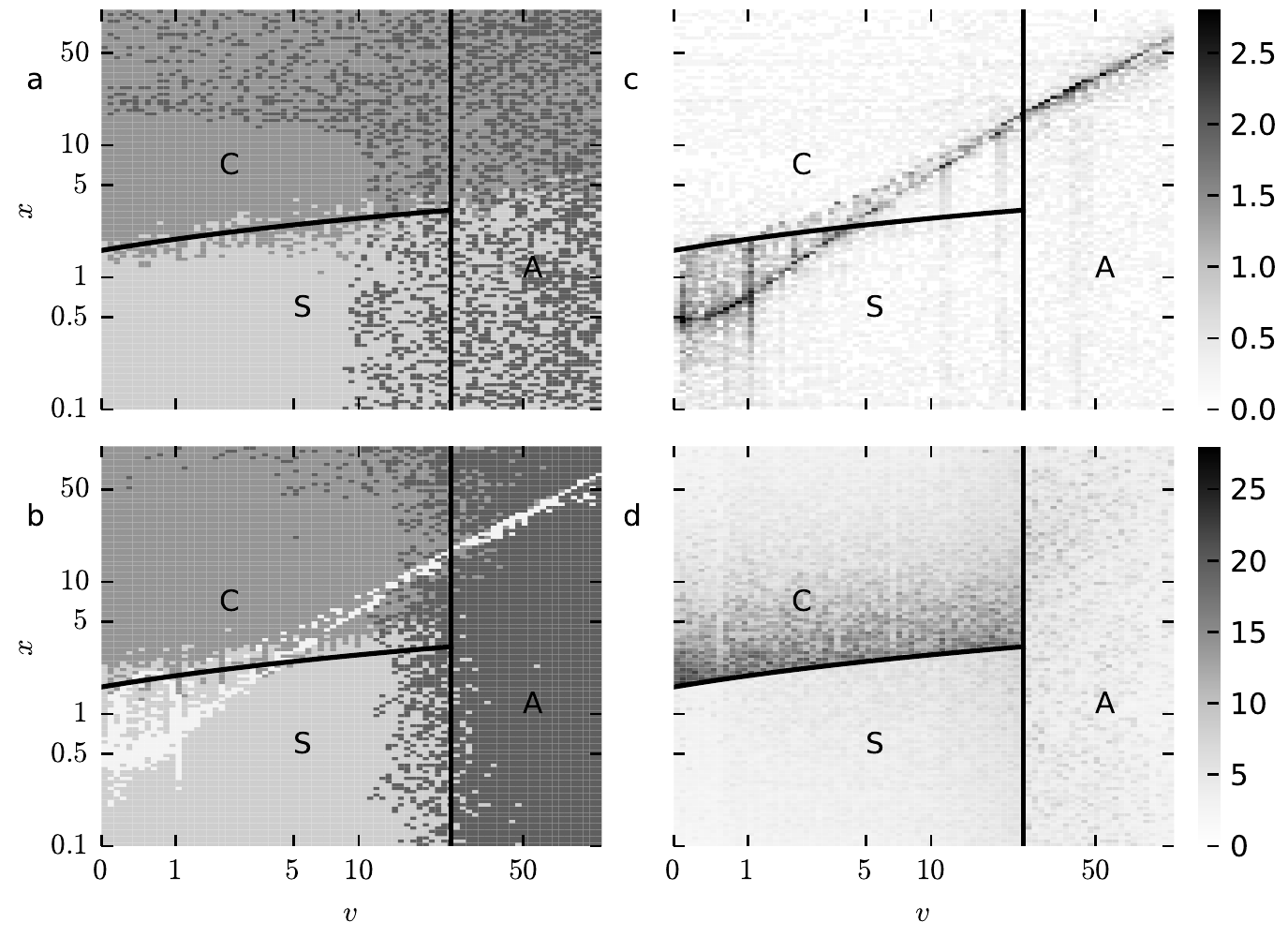}
    \caption{a) Among S, C, and A, the region with the best accuracy for S is depicted in red, the region with the best accuracy for C is in blue, and the region with the best accuracy for A is in green.
    b) Among S, C, and A, the method with the smallest computation time is shown for the condition that error is less than 1. c, d) Evaluation of theoretical error (c) and computation time (d) by combining S, C, and A.}
    \label{fig:comb}
\end{figure}

For each value of $v$ and $x$, the method with the highest accuracy among S, C, and A was identified. The results demonstrated that S and C were separated by the value of $x$. In contrast, A and S+C are distributed in the same region were sparsely and are considered to be equal in terms of accuracy (\figref{comb}a).
The criterion for separating S and C was set to $x = 1.6 + 0.5 \log \rbra{x + 1}$.

We also examined the method with the shortest computation time under the condition that the error is less than one. Accordingly, the regions of S, C, and A were clearly separated (\figref{comb}b). The criterion for separating S+C and A was $v = 25$.

When the regions were defined in such a way that the accuracy was high and the computation time was short, sufficient results were obtained for the accuracy (\figref{comb}c); however, for the computation time, there were regions near the boundary in which the computation was slow (\figref{comb}d).

\subsection{Advantage of the proposed method}
We compared the proposed method with combinations of S, C, and A.
Here, in addition to our own implementation of the combination, we also compared
the proposed methods with TensorFlowProbability (tfp), an existing implementation
of the same combination.

\begin{figure}
    \centering
    \includegraphics{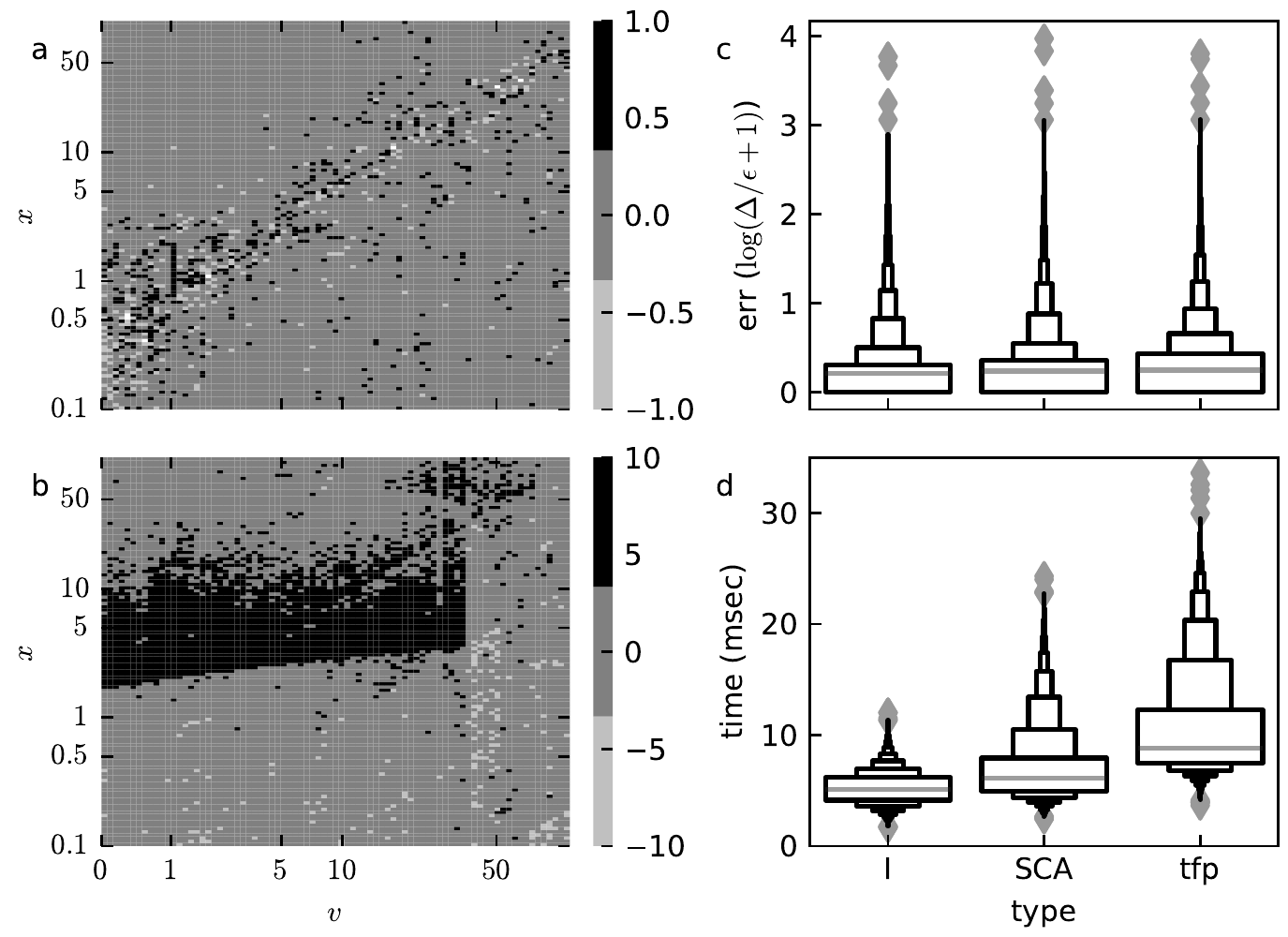}
    \caption{Comparison between the proposed method I and the SCA combination: difference in error (a) and computation time (b). Distribution of error (c) and computation time (d). tfp presents the results of TensorflowProbability.}
    \label{fig:comp}
\end{figure}

No significant difference in accuracy existed for a wide range of $v$ and $x$. The distribution of errors was not substantially different among methods (\figref{comp}a).
In addition, the distribution of the errors was not significantly different among the methods, and the mean value was the smallest for the proposed method (\figref{comp}b).

Regarding execution time, the proposed method exhibited a higher performance in a wide range (\figref{comp}c).
The computation times were stable at approximately 5 ms,
which indicates that it outperforms the other methods (\figref{comp}d).

\begin{figure}
    \centering
    \includegraphics{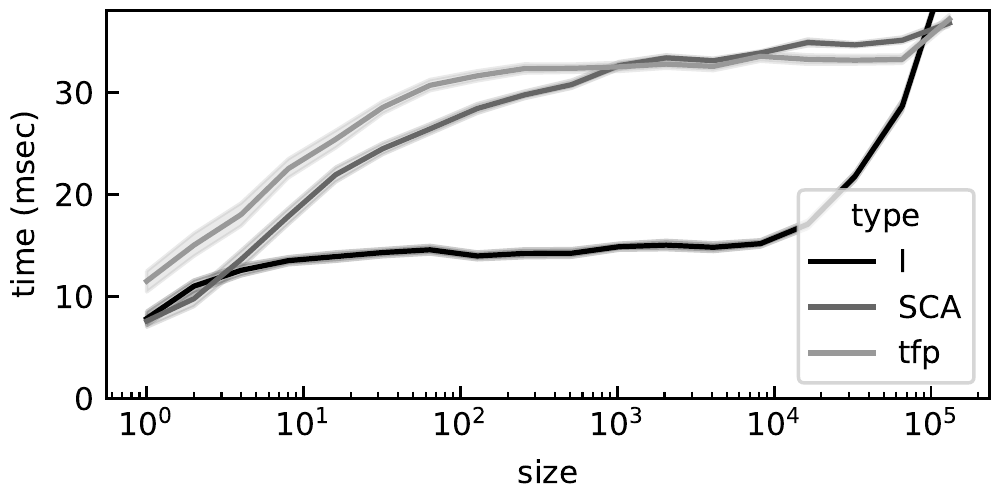}
    \caption{Computation time for parallel computing.}
    \label{fig:scale}
\end{figure}

Next, we performed vectorized computations on arrays of various sizes and measured the computation time. 
Orders $v$ and arguments $x$ were generated in the range $[0, 99]$ as $10^{2r}-1$, and $[0.1, 100]$ as $10^{3r-1}$, respectively, where $r$ denotes a uniform random number.
In this case, the proposed method could compute from  an array size of $10^1$ up to $10^4$, with almost no change in computation time  (\figref{scale}). The increase in computation time above $10^4$ can be attributed to the fact that the memory limit was reached.

For the SCA and tfp implementations, the computation time increases in the region up to $10^2$, and then becomes almost constant (\figref{scale}).
This can be expected to depend on the probability of including combinations with slow computation speed.
As aforementioned, the execution time of a vectorized computation corresponds to the worst-case computational complexity of the included computations; and 
the computation times of SCA and tfp vary for the combination of $v$ and $x$ (\figref{comp}cd).

\section{Discussion}

\subsection{Derivatives}
The derivative of the Bessel function $K_v\rbra{x}$ with respect to the argument $x$ is given by
\begin{equation}
    \frac{\partial K_{v}}{\partial x}
    = -\frac{v}{x} - \frac{K_{v-1}\rbra{x}}{K_v\rbra{x}}.\label{eq:deriv}
\end{equation}
The derivatives to $x$ of SCA and tfp are calculated using \equref{deriv}.
However, the derivative with respect to order $v$ is not provided by most libraries, and there is no known approach to obtain it indirectly.

The derivatives of the Bessel function $K_v\rbra{x}$ with respect to the order $v$ and the argument $x$
are obtained by differentiating $f_{v,x}\rbra{t}$ and integrating with $t$.
The derivatives of $f_{v,x}\rbra{t}$ are given by
\begin{gather}
    \log \rbra{-1}^m f_{v,x}^{(n,m)}\rbra{t} =
    \begin{cases}
    n \log t + m \log \cosh\rbra{t} + \log \cosh\rbra*{vt} & n~ \text{is even}\\
    n \log t + m \log \cosh\rbra{t} + \log \sinh\rbra*{vt} & n~ \text{is odd}
    \end{cases},
\end{gather}
and these also have the single range that should be integrated as the $f_{v,x}\rbra{t}$ case.
Therefore, the derivatives of $K_v\rbra{x}$ can be obtained similarly to those of $K_v\rbra{x}$.

\begin{figure}
    \centering
    \includegraphics{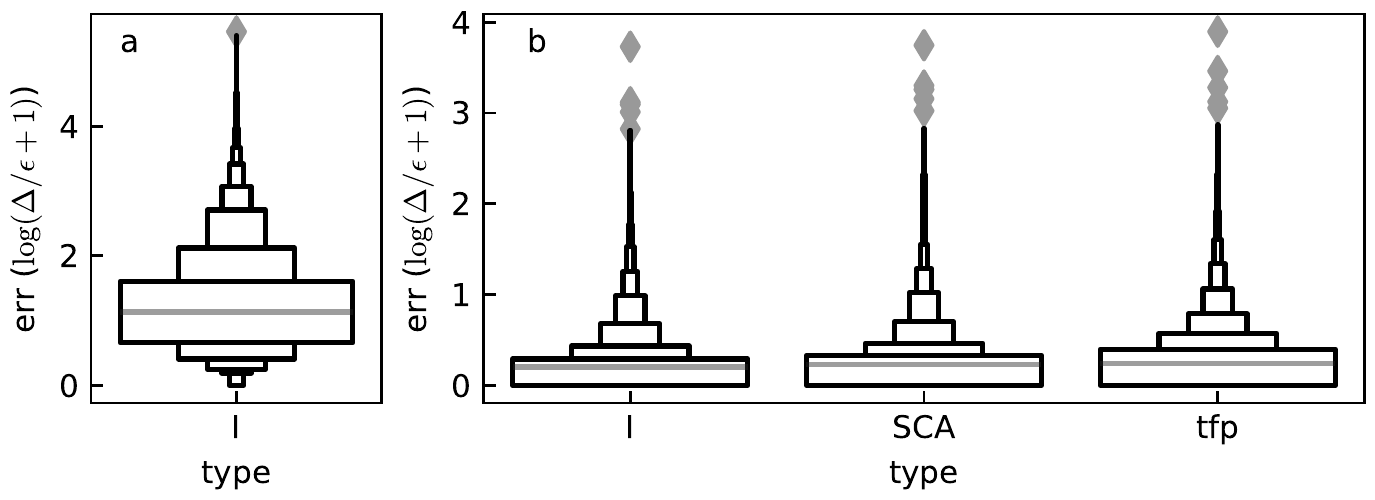}
    \caption{Differentiation errors in $v$ (a) and $x$ (b), respectively.}
    \label{fig:deriv}
\end{figure}

The derivative with respect to $x$ can be calculated with an error less than one, as $K_v\rbra{x}$ (\figref{deriv}a).
For the derivative with respect to $v$, the mean of the error is greater than one,
which indicates that the accuracy is low.
This could be attributed to the presence of the log term in the function to be integrated, and it is believed that the behavior around zero is incompletely captured.
Nevertheless, the relative error is approximately $10^14$, which is sufficiently useful,
considering that it can be calculated fast.

\subsection{Low precision floating point}

Experiments were also conducted on 32-bit low-precision floating-point systems,
and results of less than $1\epsilon$ were obtained in most regions. The computation time was almost half that of the 64-bit case, and the overall trend was also consistent. 

\section{Conclusion}

We have proposed a computational method suitable for the parallel computation of Bessel functions. The proposed method is a simple method that adds a preprocessing step to efficiently refine the integration range to the approach of computing the integral representation. Accordingly, the proposed method enables fast computation with sufficient accuracy.

This method is also applicable to the differentiation of Bessel functions, and is effective in various fields such as mathematics, physics, and statistics. In the future, it will be beneficial to investigate faster methods for the parallel processing of several related special functions.

\section*{Acknowledgements}
This work was supported by JSPS KAKENHI, Grant Number 19K12104.
We would like to thank Editage (www.editage.com) for English language editing.



\appendix

\section{Series expansion}\label{series}
When $\abs{u} \le 1/2$ and $x$ is small, $K_u\rbra{x}$ can be expanded into a series \cite{temme}:
\begin{gather}
    K_u\rbra{x} = \sum_{n=0}^\infty
    \frac{1}{n!}
    \rbra*{\frac{x}{2}}^{2n}f_n,
    \quad
    K_{u+1}\rbra{x} = 
    \sum_{n=0}^\infty
    \frac{1}{n!}
    \rbra*{\frac{x}{2}}^{2n-1}g_n,
\end{gather}
where
\begin{gather}
    f_n = \frac{\pi}{2 \sin u\pi}
    \cbra*{\rbra*{\frac{x}{2}}^{-u}\Gamma\rbra{n + 1 - u}^{-1}
    -
    \rbra*{\frac{x}{2}}^{u}\Gamma\rbra{N + 1 + u}^{-1}},\\
    g_n = p_n - n f_n,\quad
    p_n =
    \frac{1}{2\rbra*{n - u}!}
    \rbra*{\frac{x}{2}}^{-u}
    \Gamma\rbra{1+u}^{-1}.
\end{gather}
The coefficients $f_n$, $p_n$, and $q_n$ can be obtained via forward recursion:
\begin{gather}
    f_n = \frac{n f_{n-1} + p_n + q_n}{n^2 + u^2},\quad
    p_n = \frac{p_{n-1}}{n-u},\quad
    q_n = \frac{q_{n-1}}{n+u}.
\end{gather}
$K_v\rbra{x}$ = $K_{n+u}\rbra{x}$ is recursively obtained from $K_{u}$ and $K_{u+1}$:
\begin{equation}
    K_{u+1} = \frac{2 u}{x}K_{u} + K_{u-1}.\label{recu}
\end{equation}

\section{Continued fraction method}\label{continuedfraction}
Using a sequence of numbers $p_n$ containing hypergeometric functions \cite{temme, thompson}
\begin{multline}
    p_n = \rbra{2x}^n {}_2F_0\rbra*{u+n+\frac{1}{2}; u-n-\frac{1}{2}; -2x}\\
    =
    \sum_{m=0}^\infty
    \frac{\rbra{2x}^{n+m}}{m!} \prod_{i=1}^m \cbra*{\rbra*{n+i-\frac{1}{2}}^2 - u^2},
\end{multline}
the Bessel function $K_u\rbra{x}$ for a large $x$ can be defined as:
\begin{gather}
    K_u\rbra{x} = \sqrt{\frac{\pi}{2x}}e^{-x} p_0.\label{kvp}
\end{gather}
From the derivative of $p_0$ and $K_u\rbra{x}$,
an adjacent value $K_{u+1}\rbra{x}$ is also given by
\begin{gather}
    K_{u+1}\rbra{x}
    = \cbra*{\frac{1}{x}\rbra*{\frac{1}{2} + v + x} + \frac{a_1}{x}\frac{p_1}{p_0}} K_u\rbra{x}.
\end{gather}
The value of $K_v\rbra{x}$ for any order $v = n + u$ can be obtained recursively, similar to the series expansion case.

From an asymptotic expression for $p_n$
\begin{gather}
    p_{n-1} = b_n p_n - a_{n+1} p_{n+1},\label{pn3}\\
    a_n = \rbra*{n - \frac{1}{2}}^2 - u^2,\quad
    b_n = 2 \rbra{x + n},
\end{gather}
the ratio $q_n = p_n/p_0$ can be calculated recursively:
\begin{gather}
    q_{-1} = 0, \quad q_0 = 1, \quad
    q_{n+1} = \frac{q_{n-1} - b_n q_n}{a_{n+1}}.
\end{gather}
The ratios of adjacent $p_n$ are also represented using Steel's method:
\begin{gather}
    c_1 = \frac{1}{b_1},\quad
    c_{n+1} = \frac{1}{b_{n+1} - a_{n+1} c_n},\\
    r_1 = \frac{1}{b_1}, \quad
    r_{n+1} = r_n \rbra{b_{n+1} c_{n+1} - 1},\\
    \frac{p_n}{p_{n-1}} = \sum_{m=n-1}^\infty r_m.
\end{gather}

From Temme's rule 
\begin{gather}
    \sum_{m=0}^\infty d_m p_m = 1,\quad
    d_n = \frac{1}{n!} \prod_{m=1}^n a_m,\label{temme}
\end{gather}
we can calculate $p_0$ from $d_n$, $r_n$, and $p_n$:
\begin{multline}
    p_0
    = p_0 \rbra*{\sum_{n=0}^\infty d_n p_n}^{-1}
    = \rbra*{1 + \sum_{n=1}^\infty d_n \frac{p_{n-1}}{p_0} \frac{p_{n}}{p_{n-1}}}^{-1}\\
    = \rbra*{1 + \sum_{n=1}^\infty d_n q_n \sum_{m=n-1}^\infty r_m}^{-1}
    = \rbra*{1 + \sum_{m=0}^\infty r_m \sum_{n=1}^m d_n q_n}^{-1}.
\end{multline}

\section{Asymptotic expansion}\label{asymptotic}
The expansion in terms of the harmonic functions for $K_v\rbra{x}$ is given by
\cite{olver}
\begin{gather}\label{olver}
    K_v\rbra{x} = \sqrt{\frac{\pi}{2p}}\rbra*{\frac{v + p}{x}}^v e^{-p}\sum_{i=0}^\infty \sum_{j=0}^{i} c_{i,j} p^{-i} q^j ,\\
    p = \sqrt{v^2 + x^2},\quad
    q = \frac{v^2}{v^2 + x^2}.
\end{gather}
The coefficients $c_{i,j}$ can be evaluated using a recurrence relationship:
\begin{gather}\label{olver_coef}
    c_{i+1,j} =
    \rbra*{\frac{k-2}{2} + \frac{5}{8\rbra{k+1}}} c_{i,j-1}
    - \rbra*{\frac{k}{2} + \frac{1}{8\rbra{k+1}}} c_{i,j},\\
    i = 0, \dots, \quad
    j = 0, \dots, i + 1,\quad
    k = i + 2j,\\
    c_{0,-1} = c_{0,1} = 0,\quad
    c_{0,0} = 1.
\end{gather}

%% References


%\bibliographystyle{elsarticle-num}
\bibliographystyle{elsarticle-num}
\bibliography{bessel}






% \end{thebibliography}

\end{document}